\documentclass{gtpart}

\usepackage{amsmath,amsfonts,amsthm}
\usepackage{graphicx}
\usepackage{pinlabel}
\usepackage{cite}

\newtheorem{thm}{Theorem}[section]
\newtheorem{lem}[thm]{Lemma}
\newtheorem{prop}[thm]{Proposition}
\newtheorem{cor}[thm]{Corollary}

\newtheorem{constr}[thm]{Construction}

\theoremstyle{definition}
\newtheorem{defn}[thm]{Definition}

\theoremstyle{remark}
\newtheorem*{rmk}{Remark}

\newenvironment{pf}{ \begin{proof} }{ \end{proof} }

\DeclareMathAlphabet\EuScript{U}{eus}{m}{n}
\SetMathAlphabet\EuScript{bold}{U}{eus}{b}{n}

\DeclareFontFamily{U}{eus}{\skewchar\font'60}%
\DeclareFontShape{U}{eus}{m}{n}{<-6>eusm5<6-8>eusm7<8->eusm10}{}%
\DeclareFontShape{U}{eus}{b}{n}{<-6>eusb5<6-8>eusb7<8->eusb10}{}%

\DeclareMathOperator{\im}{im}
\DeclareMathOperator{\rank}{\mathrm{rank}}
\DeclareMathOperator{\sym}{Sym}
\DeclareMathOperator{\hilb}{Hilb}

\DeclareMathOperator{\spinc}{\mathrm{Spin}^c}
\DeclareMathOperator{\interior}{int}
\DeclareMathOperator{\crit}{crit}
\newcommand{\id}{\mathrm{id}}

\begin{document}
\title[Hamiltonian handleslides for Heegard Floer homology]{Hamiltonian handleslides for Heegaard Floer homology}
\author{Timothy Perutz}
\date{December 19, 2007}
\address{Department of Mathematics, Columbia University, 2990 Broadway, New York, NY 10027, USA.}
\email{perutz@math.columbia.edu}

\begin{abstract}
A $g$-tuple of disjoint, linearly independent circles in a Riemann surface $\Sigma$ of genus $g$ determines a `Heegaard torus' in its $g$-fold symmetric product. Changing the circles by a handleslide produces a new torus. It is proved that, for symplectic forms with certain properties, these two tori are Hamiltonian-isotopic Lagrangian submanifolds. This provides an alternative route to the handleslide-invariance of Ozsv\'ath--Szab\'o's Heegaard Floer homology.
\end{abstract}
\maketitle
\section{Introduction}
\subsection{Handlebodies and handleslides}
Let $\Sigma$ be a surface of  genus $g\geq 1$. We can express $\Sigma$ as the boundary $\partial U$ of a 3-dimensional handlebody $U$ by choosing $g$ disjoint, embedded circles, $(\gamma_1,\dots,\gamma_g)$, linearly independent in homology: each of these is filled so as to bound a disc in $U$. Conversely, the handlebody $U$ determines the equivalence class of the $g$-tuple of attaching circles $(\gamma_1,\dots,\gamma_g)$ under the equivalence relation generated by isotopies, permutations and \emph{handleslides} (see, for instance, \cite{GS}).  

Handleslides are not possible when $g=1$ (the attaching circle for a genus 1 handlebody is unique up to isotopy), so from now on we shall assume that $g\geq 2$. A handleside is a replacement
\[ (\gamma_1,\gamma_2,\dots,\gamma_g)\leadsto (\gamma_0,\gamma_2,\dots, \gamma_g),\] 
where $\gamma_0$ is disjoint from $\gamma_2\cup \dots \cup \gamma_g$ and isotopic to another circle $\gamma_0'$, disjoint from $\gamma_1\cup \gamma_2$, such that $ \gamma_0'\cup \gamma_1\cup \gamma_2$ bounds an embedded pair of pants in $\Sigma$; see Figure 1. This condition implies that the circles $(\gamma_0,\gamma_2,\dots, \gamma_g)$ are again linearly independent in homology, and hence determine a handlebody $U'$. 
There is a diffeomorphism $U\to U'$ acting as the identity on $\partial U = \Sigma = \partial U'$, since, once the $g-1$ curves $(\gamma_2,\dots,\gamma_g)$ have been collapsed,  $\gamma_0$ becomes isotopic to $\gamma_1$. 

\begin{figure}[t!]
\labellist
\small\hair 2pt
\pinlabel $\gamma_2$ at 220 325
\pinlabel $\gamma_2$ at 220 115
\pinlabel $\gamma_0'$ at 165 245
\pinlabel $\gamma_0$ at 165 32
\pinlabel $\gamma_1$ at 330 245
\pinlabel $\gamma_1$ at 323 50
\pinlabel $A$ at 220 60
\pinlabel $D$ at 350 20
\endlabellist
\centering
\includegraphics[width=0.7\textwidth]{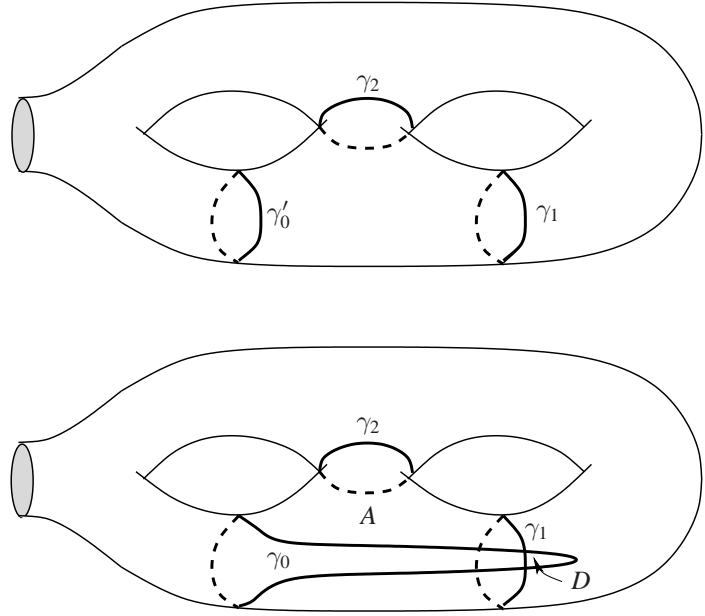}
\caption{A handleslide in a punctured genus 2 surface. In the upper diagram, the three curves $\gamma_0'$, $\gamma_1$ and $\gamma_2$ bound an embedded pair of pants. In the lower diagram, $\gamma_0$ intersects $\gamma_1$ in two points; this pattern of intersection is the one we shall use in the proof of our theorem.}
\end{figure}
\subsection{Heegaard tori}
Fix, once and for all, a complex structure $j$ on $\Sigma$. The $g$-fold symmetric product $\sym^g(\Sigma)$ is then a complex manifold and hence also a differentiable manifold. In working with $\sym^g(\Sigma)$ one should be aware that its differentiable structure depends on $j$ but its diffeomorphism type does not. The $g$-tuple of disjoint attaching circles $(\gamma_1,\dots,\gamma_g)$ determines an embedded torus in $\sym^g(\Sigma)$,
\[ T_1 = \pi( \gamma_1 \times \gamma_2 \times\dots \times \gamma_g) , \]
its \emph{Heegaard torus}. Here $\pi\colon \Sigma^{\times g}\to \sym^g(\Sigma)$ is the quotient map; it restricts to an embedding on $ \gamma_1 \times \gamma_2  \times \dots \times \gamma_g$ since the circles $\gamma_i$ are disjoint. Likewise, the $g$-tuple 
$(\gamma_0,\gamma_2,\dots,  \gamma_g) $ determines a torus
\[T_0 =  \pi( \gamma_0 \times \gamma_2  \times \dots \times \gamma_g).\]
These tori are totally real, that is, their tangent spaces contain no complex lines; 
indeed, if we choose an area form $\alpha $ on $\Sigma$ then the complement of the diagonal in $\sym^g(\Sigma)$ has an induced symplectic form $\pi_*(\alpha^{\times g})$, and this form makes the two tori Lagrangian. Note, however, that the push-forward of a 2-form by a finite holomorphic map is not globally a smooth 2-form but rather a current, with singularities along the branch locus---in this case, the diagonal. 

We now recall that for $n>1$ one has $H^2(\sym^n(\Sigma);\Z) \cong H^0(\Sigma;\R) \oplus \Lambda^2 H^1(\Sigma;\Z)$ by an isomorphism equivariant under the actions of the mapping class group on the two sides. Write $\eta$ for the class corresponding to $1\in H^0(\Sigma;\Z)$, or for its image in $H^2(\sym^n(\Sigma);\R)$ (this leaves a sign ambiguity which can be fixed by giving an alternative description: $\eta$ is Poincar\'e dual to the divisor $\sym^{g-1}(\Sigma)$, embedded by $D\mapsto p + D$ for some fixed $p$). When $g>0$, write $\theta$ for  (the image in real cohomology of) the generator of the cyclic group of classes in the summand $\Lambda^2 H^1(\Sigma;\Z)$  invariant under the mapping class group. This time the sign ambiguity can be fixed by saying that $\theta$ is the pullback by the Abel--Jacobi map of an ample class (the theta-divisor) on the Jacobian. Suitable references include \cite{ACGH, BT}.

Our first, preliminary result is as follows.
\begin{prop}\label{forms}
Suppose $(\gamma_0,\gamma_2,\dots,\gamma_g)$ is a $g$-tuple of disjoint, linearly independent attaching circles, and $(\gamma_1,\gamma_2,\dots, \gamma_g)$ another $g$-tuple differing from the first by a handleslide. For real numbers $\lambda$ with sufficiently small absolute value, and arbitrary positive area forms $\alpha$ on $\Sigma$, there exist symplectic forms $\omega_\lambda$ on $\sym^g(\Sigma)$, taming the complex structure $\sym^g(j)$, such that
\begin{enumerate}
\item $\omega_\lambda$ agrees with the product form $\pi_*(\alpha^{\times g})$ in neighbourhoods of the Heegaard tori $T_0$ and $T_1$; and
\item
$\omega_\lambda$ represents the class $\eta+\lambda \theta$.
\end{enumerate}
\end{prop}
One can in fact take $\omega_\lambda$ to be K\"ahler. A short proof can be given using smoothing theory for currents, as in the author's earlier unpublished note \cite{Pe1} which has now been incorporated as the last section (\ref{smoothing}) of this paper. An immediate consequence of this smoothing theory (Corollary \ref{smoothing cor}) is that there are forms satisfying the proposition in the class $\eta$. From this it is easy to deduce the result: one can perturb such a form by adding a small multiple of a closed $(1,1)$-form supported near the diagonal $\delta$ and Poincar\'e dual to $\delta$. Since $PD [\delta] = \eta  - \theta$,\footnote{The tangent bundle to $\sym^n(\Sigma)$ can be described algebro-geometrically as $\mathcal{O}_\delta(\delta)$. Hence $PD [\delta]= c_1(T \sym^n(\Sigma))=(n+1-g)\eta  - \theta$ \cite{ACGH}.} these perturbed forms do the job.

This paper actually contains a second proof of the proposition, independent from that given in Section \ref{smoothing}. The second proof is, however, valid only when $\lambda>0$, and the forms it produces are not necessarily K\"ahler. It can be found in Section \ref{isotopies}, subsumed in the proof of our main theorem, which we now proceed to state.

\subsection{The main theorem} 
For a form $\omega_\lambda$ as in Proposition \ref{forms}, the tori $T_0$ and $T_1$ are Lagrangian. One may ask whether they are isotopic through Lagrangians, or indeed Hamiltonian-isotopic. To make these into well-posed questions, we suppose that $\gamma_0$ and $\gamma_1$ intersect transversely in precisely two points, as in Figure 1. We also fix $\alpha$. The answers to the questions are then independent of the particular form $\omega_\lambda$ and depend only on the parameter $\lambda$. Indeed, it follows from  Moser's lemma and the convexity of the set of symplectic forms taming $\sym^g(j)$ that, when such forms are cohomologous, they are related by symplectomorphisms supported away from $T_0\cup T_1$. 

We answer these questions when $\lambda>0$. The answers depend on the areas $\int_D{\alpha}$ and $\int_A{\alpha}$ of the half-disc $D$ enclosed by arcs in $\gamma_0$ and $\gamma_1$ and the the annular region $A$ enclosed by $\gamma_2$ and the remaining portions of $\gamma_0$ and $\gamma_1$ (see Figure 1).

\begin{thm}\label{main}
Suppose $(\gamma_0,\gamma_2,\dots,\gamma_g)$ is a $g$-tuple of disjoint, linearly independent attaching circles, and $(\gamma_1,\gamma_2,\dots, \gamma_g)$ another $g$-tuple differing from the first by a handleslide. Assume that $\gamma_0\cap \gamma_1$ is a transverse intersection consisting of precisely two points. Fix a positive area form $\alpha$ on $\Sigma$, and let $\omega_\lambda$ be a K\"ahler  form satisfying the conclusions of the last proposition. Suppose that $\lambda$ is strictly positive. Then
\begin{enumerate}
\item
$T_0$ is Lagrangian-isotopic to $T_1$.
\item
$T_0$ is Hamiltonian-isotopic to $T_1$ if and only if $\int_D{\alpha} = (1+\lambda) \int_{A}{\alpha} - 2 \lambda.$
\item
A Lagrangian isotopy (or indeed a Hamiltonian one, when the area constraint is satisfied) may be constructed as the product of the constant isotopy on the common factor $\gamma_3\times\dots \times \gamma_g$ of the two tori and an isotopy in $\sym^2(\Pi)$, where $\Pi\subset \Sigma$ is an (arbitrary) embedded pair of pants containing $\gamma_0 \cup \gamma_1 \cup\gamma_2$ in its interior.
\end{enumerate}
\end{thm}

The proof of Lagrangian isotopy is constructive. The idea is to exploit the fact that $\gamma_0$ and $\gamma_1$ become isotopic in $\Sigma'$, the surface obtained by surgering out $\gamma_2$. The genus 2 case contains the heart of the matter. When $\Sigma$ has genus 2, the surgery is mimicked symplectically by associating with $\gamma_2$ a hypersurface $V\subset \sym^2(\Sigma)$ which is the total space of an $S^1$-bundle $\rho\colon V\to \Sigma'$. The pullback by $\rho$ of an area form on $\Sigma'$ agrees with the symplectic form on $\sym^2(\Sigma)$ restricted to $V$. The preimage under $\rho$ of a Lagrangian isotopy between the images of $\gamma_0$ and $\gamma_1$ in $\Sigma'$ is then a Lagrangian isotopy in $\sym^2(\Sigma)$. One needs to be able to exert enough  control over $V$ to arrange that the isotopy begins at $\gamma_0\times \gamma_2$ and ends at $\gamma_1\times \gamma_2$. 

Note that the forms used in the proof are \emph{not} the perturbed, smoothed currents mentioned above, but rather, symplectic structures arising from presentations of symmetric products as symplectic sums.  However, as already noted, the existence of Lagrangian or Hamiltonian isotopies is a function of $\alpha$ and $\lambda$ and not of the particular form $\omega_\lambda$.

\begin{rmk} 
It is an open question whether one can find a form making $T_0$ and $T_1$ into Hamiltonian-isotopic Lagrangians when $\lambda=0$ and $\int_D{\alpha} = \int_{A}{\alpha}$. It is conceivable that we have here an instance of \emph{symplectic fragility}, the phenomenon whereby non-isotopic Lagrangians become isotopic as soon as a perturbation parameter is switched on \cite{Sei}, though this seems rather unlikely since we shall find that the Floer-theoretic properties of $(T_0,T_1)$ do not change as $\lambda$ decreases to $0$. In the gauge-theoretic interpretation of symmetric products as moduli spaces of vortices, the natural symplectic form represents (up to scale) a class $\eta+\lambda \theta$ where $\lambda$ decreases to $0$ as the `stability parameter' tends to infinity \cite{Pe3}.
\end{rmk}

We shall use Theorem \ref{main} to re-prove the handleslide-invariance of 
Ozsv\'ath--Szab\'o's Heegaard Floer homology \cite[Section 9]{OS} (the argument still draws on finiteness lemmas from \cite{OS}). The precise statement uses notions from \cite{OS} which we do not recapitulate here.
\begin{cor}\label{HeegCor}
Let $(\Sigma, \boldsymbol{\alpha}, \boldsymbol{\beta},z)$ be a pointed Heegaard diagram,  weakly admissible for a $\spinc$-structure $\mathfrak{s}$, and suppose that $\boldsymbol{\alpha}'$ differs from $\boldsymbol{\alpha}$ by a handleslide. Suppose that the circle that is being replaced intersects its replacement non-trivially in the pattern of Theorem \ref{main}.  Assume that $z$ lies outside the handlesliding region.  Then there is a  canonical isomorphism $HF^+(\boldsymbol{\alpha},  \boldsymbol{\beta}, \mathfrak{s}) \to HF^+(\boldsymbol{\alpha'},  \boldsymbol{\beta},\mathfrak{s})$ between the Heegaard Floer homology groups for the $\spinc$-structure $\mathfrak{s}$, and similarly for $HF^\infty$, $HF^-$ and $\widehat{HF}$.
\end{cor}
Referring to Figure 1, the `handlesliding region' means the union of the region bounded by $\gamma_2$ and arcs in $\gamma_0$ and $\gamma_1$ and that bounded by arcs in $\gamma_0$ and $\gamma_1$.

\emph{Plan of the paper.} In Section \ref{review} we review the relationship between degenerations of symplectic manifolds and the symplectic sum operation; this is then applied to symmetric products (Section \ref{symm}), and employed to construct the Lagrangian and Hamiltonian isotopies of the theorem (Section \ref{isotopies}). Section \ref{non-iso} proves the Hamiltonian non-isotopy clause, and Section \ref{Heegaard} covers the application to Heegaard Floer homology. 
The symplectic forms used up to that point tame the natural complex structure on $\sym^g(\Sigma)$ but are not actually K\"ahler. It has already been mentioned that in Section \ref{smoothing} we give a different construction, by smoothing currents, which does produce K\"ahler forms and has some independent interest.

\subsection*{Acknowledgements}
The question of whether handleslides lead to Hamiltonian-isotopic Heegaard tori was posed to me by Peter Ozsv\'ath. Thanks to him, and to Ron Fintushel, Robert Lipshitz and Ivan Smith for productive conversations. I am grateful to Selman Akbulut and the other organisers of the 14th G\"okova Geometry--Topology Conference for a mathematically stimulating week in an idyllic setting.

\section{Degeneration and symplectic sums}\label{review}
In this section we review the notion of symplectic sum and its relation to degenerations of symplectic manifolds \cite{Gom, IP}. So as to simplify the definitions (slightly) we restrict the discussion to complex manifolds. 
\begin{defn}
If $M$ is a complex $n$-manifold, a \emph{nodal degeneration} of $M$ is a pair $(E,\pi)$ where $E$ is a complex $(n+1)$ manifold and $\pi\colon E\to \C$ a holomorphic map which is a topological fibre bundle over $\C^*$, with fibre $E_1=M$ (we write $E_t$ for $\pi^{-1}(t)$). Thus the critical locus $C=\crit(\pi)$ is contained in the zero-fibre $E_0$. We suppose that this locus is smooth and $(n-1)$-dimensional, and that the (complex) Hessian on $N_{C/E}$ is everywhere non-degenerate.
\end{defn}
As a complex space, $E_0$ has a normal crossing along $C$: locally near a  point of $C$,  $E_0$ is equivalent to a neighbourhood of the origin in $\C^{n-1} \times \{ z_1 z_2 =0  \}\subset \C^{n-1}\times\C^2$. Let $n\colon \tilde{E}_0\to E_0$ be the normalisation of $E_0$. The normalisation is an intrinsic construction but, put simply, it is the complex manifold obtained by replacing each patch $\C^{n-1} \times \{ z_1 z_2 =0  \}$ by $\C^{n-1}\times (\C \amalg \C)$. The preimage $n^{-1}(C)\subset \tilde{E}_0$ is the disjoint union of two divisors $Z_1$ and $Z_2$ which are identified with one another via $n$. Their normal bundles are dually paired with one another by means of the Hessian $D^2\pi$ on $N_{C/E}$. Choose  closed tubular neighbourhoods $U_1$ and $U_2$ in $\tilde{E}_0$ for $Z_1$ and $Z_2$ respectively, and put $V_i = \partial U_i$, so that $V_i$ is an $S^1$-bundle over $Z_i$. There is an orientation-reversing diffeomorphism $\sigma \colon V_1 \to V_2$, which arises because because we can think of $V_1$ as the equator in the $\mathbb{CP}^1$-bundle $\mathbb{P}(N_1\oplus \mathcal{O})\to Z_1=C$, and $V_2$ as the equator in $\mathbb{P}(\mathcal{O} \oplus N_2^*)\to Z_2= C$; but using the Hessian we have
\[ \mathbb{P}(N_1\oplus \mathcal{O})=\mathbb{P}( (N_1\otimes N_2^*) \oplus N_2^*) \cong  \mathbb{P}( \mathcal{O}\oplus N_2^*). \]

This discussion makes sense symplectically. If we are given a symplectic form $\Omega$ on $E$, taming the complex structure, the manifold $\tilde{E}_0$ inherits a symplectic form $\omega$ for which $n_*(\omega|_{Z_1}) =n_*(\omega|_{Z_2})$. We are thus in a position to form Gompf's \emph{symplectic sum}  
\[  \#_{Z_1\sim Z_2} \tilde{E}_0\] 
as in \cite{Gom}, the symplectic manifold obtained from $\tilde{E}_0 \setminus (\interior{U}_1 \cup \interior{U}_2)$ by gluing $V_1$ to $V_2$ via $\sigma$. 

\begin{prop}\label{sum}
The fibre $M=E_1$ of  a nodal degeneration $(E,\pi)$, when equipped with the restriction of a K\"ahler form $\Omega$ on $E$, is symplectomorphic to the symplectic sum $  \#_{Z_1\sim Z_2} \tilde{E}_0$.
\end{prop}
\begin{pf}[Sketch of proof]
We refer to \cite{IP} for a complete treatment; here we aim mainly to draw attention to an important aspect of the geometry of degenerations, their \emph{vanishing cycles}. 

The vanishing cycle $V\subset M$ is the set of points $x$ such that symplectic parallel transport $M = E_1 \to E_t$ over the ray $[t,1]\subset \C$ lands in the critical set $C$ in the limit $t\to 0^+$. The limiting parallel transport defines an $S^1$-bundle $\rho\colon V\to C$, and the crucial point for us will be that 
\[ \Omega|_V = \rho^*(\Omega|_C).\] 
Correspondingly, in the symplectic sum $\#_{Z_1\sim Z_2}  \tilde{E}_0 $, let $V'$ denote the common image of $V_1$ and $V_2$. Since it arises as the boundary of a tubular neighbourhood of $C$, it comes equipped with an $S^1$-bundle $\rho' \colon V'\to C$. 

Notice that $ (\#_{Z_1\sim Z_2}  \tilde{E}_0 ) \setminus V'$ is naturally identified symplectically with $E_0 \setminus C$. On the other hand, symplectic parallel transport into $E_0$ defines a symplectomorphism $M\setminus V \to E_0 \setminus C$. Combining these observations we find a symplectomorphism $M\setminus V \cong (\#_{Z_1\sim Z_2}  \tilde{E}_0 ) \setminus V'$. But both $V$ and $V'$ are coisotropic submanifolds in symplectic manifolds (in fact, they are $S^1$-bundles over a common symplectic manifold). There is a diffeomorphism $V \to V'$ covering the identity map on $C$, and by the coisotropic neighbourhood theorem, this extends to a symplectomorphism between neighbourhoods of $V$ and $V'$. It remains to see that these can be made compatible with  the identification $M\setminus V \cong (\#_{Z_1\sim Z_2}  \tilde{E}_0 ) \setminus V'$. This requires closer inspection of the local model for the symplectic sum, and we shall not give the details here (see \cite{IP}; also compare \cite[Section 2]{Pe2}).
\end{pf}


\section{Symmetric products as symplectic sums}\label{symm}

Let $\gamma \subset \Sigma$ be a simple closed curve, and let $\Sigma_\gamma$ be the surface obtained by excising a tubular neighbourhood of $\gamma$ and gluing in a pair of discs. Let $p$ and $q$ be points in the respective interiors of those discs.
Fix complex structures on $\Sigma$ and $\Sigma_\gamma$, and an integer $n\geq 1$. 
\begin{constr}
Define two maps
\[  i_p, i_q \colon \sym^{n-1}(\Sigma_\gamma)  \to \sym^n(\Sigma_\gamma) \]
by $i_p(D)=p+D$, $i_q(D)=q+ D$.  Denote by $\mathbf{S}_n(\Sigma_\gamma)$ 
the complex  blow-up of $\sym^n(\Sigma_\gamma)$ along the locus 
$p+q + \sym^{n-2}(\Sigma_\gamma)$. Let
\[  \tilde{\imath}_p, \tilde{\imath}_q \colon \sym^{n-1}(\Sigma_\gamma)  \to \mathbf{S}_n(\Sigma_\gamma)  \] 
be the natural lifts of $i_p$ and $i_q$. The images $Z_p=\im \tilde{\imath}_p$ and $Z_q=\im \tilde{\imath}_q$ are disjoint.  Let $\tau =\tilde{\imath}_p \circ \tilde{\imath}_q^{-1} \colon Z_q\to Z_p$.  There is a natural isomorphism $\tau^* N_{Z_p/\mathbf{S}_n(\Sigma_\gamma)} \cong N_{Z_q/\mathbf{S}_n(\Sigma_\gamma)}^*$ (see \cite[Section 3]{Pe2}). Thus $\tau$ identifies the divisors $Z_p$ and $Z_q$, and under this identification, their normal bundles are dual. We may therefore form the symplectic sum $ \#_{Z_p\sim Z_q} \mathbf{S}_n(\Sigma_\gamma) $ (as a smooth manifold; we will say more about symplectic forms presently).
\end{constr}

\begin{prop}
There is a diffeomorphism
\[  \phi \colon  \#_{Z_p\sim Z_q} \mathbf{S}_n(\Sigma_\gamma) \cong \sym^n(\Sigma),   \]
canonical up to isotopy.
\end{prop}

\begin{pf}
The diffeomorphism-type of a symmetric product of a Riemann surface does not depend on its complex structure. Moreover, because the space of complex structures is simply connected, there are canonical diffeomorphisms between them (up to isotopy).  Thus we can choose complex structures on $\Sigma$ and $\Sigma_\gamma$ as we wish.

To prove the proposition, we invoke Proposition \ref{sum}. What we are asserting is that there is a nodal degeneration of $\sym^n(\Sigma)$ to a complex space whose normalisation is $\mathbf{S}_n(\Sigma_\gamma)$, in which $Z_p$ and $Z_q$ are the two preimages of the normal crossing divisor, and that the pairing of normal bundles is induced by the Hessian along the critical locus of the degeneration.

Let $E\to \C$ be a nodal degeneration of $\Sigma$, i.e. a holomorphic Lefschetz fibration such that $E_1=\Sigma$, with a single critical point lying over $0$. We arrange that the vanishing cycle, taken along the vanishing path $[0,1]$, is (the isotopy class of) $\gamma$. One can then form an associated family, $\hilb^n(\pi)\to \C$, the relative Hilbert scheme of $n$ points, as in Donaldson--Smith \cite{DS}. Their crucial observation (later re-proved by Ran \cite{Ran}) is that this space is globally smooth; indeed, it is a nodal degeneration of its smooth fibre. 
There is a natural morphism $\hilb^n(\pi)\to \sym^n(\pi)$ to the relative symmetric product, whose restriction to $\C^*$ is biholomorphic onto the restriction of $\sym^n(\pi)$ to  $\C^*$; hence $\hilb^n(\pi)$ is a nodal degeneration of $\sym^n(\Sigma)$. 

The structure of the relative Hilbert scheme was described in elementary terms in \cite{Pe2}; in particular, it was explained that the normalisation of the zero-fibre $\hilb^n(E_0)$ is precisely $\mathbf{S}_n(\widetilde{E}_0)$ where $\widetilde{E}_0\to E_0$ is the normalisation of the zero-fibre of $E\to \C$. The two distinguished divisors inside it are the divisors $Z_p$ and $Z_q$ described in the discussion above (with $p$ and $q$ the points of $\widetilde{E}_0$ lying over the node of $E_0$), so $\hilb^n(E_0)$ is obtained by identifying $Z_p$ with $Z_q$ via $\tau$.
Thus our result follows from the general theory of Section \ref{review}.
\end{pf}
A case to keep in mind is the symmetric square, $\sym^2(\Sigma)$. In this case, $\mathbf{S}_2(\Sigma_\gamma)$ is the blow-up of $\sym^2(\Sigma_\gamma)$ at the point $p+q$. The two divisors $Z_p$ and $Z_q$ are disjoint copies of 
$\Sigma_\gamma$. If $\gamma$ is a separating curve, so that $\Sigma_\gamma $ is a disjoint union $\Sigma_1 \amalg \Sigma_2$, the assertion is that $\sym^2(\Sigma)$ is the fibre sum
\[  \sym^2(\Sigma_1)  \#_{\Sigma_1}  
(\widetilde{\Sigma_1 \times \Sigma_2})  
\#_{ \Sigma_2}   \sym^2(\Sigma_2) .  \]
Here $\widetilde{\Sigma_1 \times \Sigma_2}$ is the blow-up of $\Sigma_1 \times \Sigma_2$ at $(p,q)$; the proper transforms of the factors in the product give embeddings of $\Sigma_1$ and $\Sigma_2$ into the blow-up, each of self-intersection $-1$. On the other hand, $\sym^2(\Sigma_1)$ contains a copy of $\Sigma_1$ of self-intersection $+1$, embedded by the map $x\mapsto p+x$. Likewise, $\sym^2(\Sigma_2)$ contains a copy of $\Sigma_2$ of self-intersection $+1$, embedded by the map $x\mapsto q+x$. The fibre sums $\#_{\Sigma_1}$ and $\#_{\Sigma_2}$ are taken along these two pairs of surfaces.
(Note: The manifold appearing here is written as a sum of three different manifolds along two pairs of surfaces, whereas the proposition expresses it as a self-fibre sum of a disconnected manifold along a disconnected surface.)

\begin{rmk}
This remark is due to R. Fintushel.\footnote{The conclusion was also known to I. Baykur. Any mistakes in the argument are the author's responsibility.}
When $\Sigma_2$ has genus 1, the description can be simplified because the second fibre sum has no topological effect. Thus, if we write $\Sigma=\Sigma' \# T^2$ (connected sum) we have
\[ \sym^2(\Sigma) \cong \sym^2(\Sigma') \#_{\Sigma'} \widetilde{( \Sigma' \times T^2)}.  \]  
Indeed, we know that $\sym^2(\Sigma) \cong \sym^2(\Sigma') \#_{\Sigma'} \widetilde{( \Sigma' \times T^2)} \# \sym^2(T^2)$. The summand $\sym^2(T^2)$ is the total space of a non-trivial $S^2$-bundle over $T^2$ (the bundle projection is the Abel--Jacobi map). To perform the second fibre sum, we remove the tubular neighbourhood of a square $-1$ torus in $\sym^2(\Sigma)' \# (\widetilde{\Sigma' \times T^2})$ (the proper transform of $T^2 \times\{\mathrm{pt}.\}$ in the second summand). We glue in the complement of a tubular neighbourhood of a square $+1$ section $s_0$ of $\sym^2(T^2)\to T^2$. This complement is a tubular neighbourhood of another section $s_\infty$, which must have square $-1$ (since the fact that $[s_\infty] \cdot [s_0] = 0$ forces $[s_\infty]=[s_0]-[\mathrm{fibre}]$). Thus we glue into $\sym^2(\Sigma)' \# (\widetilde{\Sigma' \times T^2})$  the same piece that we removed, without changing the gluing. 
\end{rmk}

We now describe the how the symplectic sum description affects the relevant cohomology classes, adding subscripts to the notation to track which surface we are considering.

\begin{lem}\label{coh}
Assume $\Sigma_\gamma$ is connected. Let $c_\lambda \in H^2( \mathbf{S}_n(\Sigma_\gamma);\R)$ be the pullback of  $\eta_{\Sigma_\gamma}+\lambda\theta_{\Sigma_\gamma}$ from $\sym^n(\Sigma_\gamma)$ to its blow-up, minus $\lambda$ times the class Poincar\'e dual to the exceptional divisor. Under the diffeomorphism $\phi$ of the previous proposition, $c_\lambda$ pulls back to $\eta_\Sigma + \lambda\theta_\Sigma$.
\end{lem}
Closely related results were obtained in \cite{Pe2}, notably Proposition 3.14. 

The relevance of $c_\lambda$ is that it will be the class of our symplectic form. See \cite{MS,Voi}  for accounts of blowing up symplectic or K\"ahler manifolds.
\begin{pf}
We shall give a direct proof in a `prototypical' case, then explain how to obtain the general result from this.

Consider $\sym^2(T^2)$, thought of as $\#_{\mathbb{P}^1\sim \mathbb{P}^1} \widetilde{\mathbb{P}}^2$ (symplectic sum along two square-zero $2$-spheres). One has $H_2(\widetilde{\mathbb{P}}^2;\Z)=\Z^2$, the generators being the classes $e$ of the  exceptional divisor and $s$ of the proper transform of a (generic) line in $\mathbb{P}^2$. The fibre sum description corresponds to a nodal degeneration $\hilb^2_\C(E)\to \C$; the homology cycles in $\widetilde{\mathbb{P}}^2$ project to cycles in the central fibre of $\hilb^2_\C(E)$, and (by thinking of the symplectic sum construction) one can see explicitly that these projected cycles are homologous in $\hilb^2_\C(E)$ to cycles in $\sym^2(T^2)$. 

Indeed, $e$ intersects each of the two square-zero $\mathbb{P}^1$s in a point, and so deleting a neighbourhood of  $\mathbb{P}^1\amalg \mathbb{P}^1$ has the effect of removing two discs from $e$; in the symplectic sum, the two boundary circles of this 2-holed sphere are identified to form a torus $\tilde{e}$ in $\#_{\mathbb{P}^1\sim \mathbb{P}^1} \widetilde{\mathbb{P}}^2 = \sym^2(T^2)$, of square $-1$, homologous in  $\hilb^2_\C(E)$ to the projection of $e$.
Starting from $s$ one similarly gets a  square $+1$ torus in $\sym^2(T^2)$. Note that both $\tilde{e}$ and $\tilde{s}$ can be thought of as sections of the Abel--Jacobi fibration $\sym^2(T^2)\to T^2$, $[x,y]\mapsto x+y$ ($x+y$ means the sum under the group law of an elliptic curve). One now sees that $\tilde{e}$ is dual to $\theta_\Sigma$, $\tilde{s}$ to $\eta_\Sigma$. Thus in this case, the symplectic class $(\eta_{\Sigma_\gamma} +\theta_{\Sigma_\gamma} )+ \lambda e $  does indeed correspond to $\eta_\Sigma+\lambda\theta_\Sigma$ (but here $\theta_{\Sigma_\gamma}=0$).

To obtain the general formula, it suffices to do so when $n=2$. Indeed, the construction of the degeneration of $\Sigma$, and hence of its relative Hilbert scheme, is natural under diffeomorphisms of $\Sigma$ which act trivially near $\gamma$. Thus $\phi^*c_\lambda$ has to be invariant under the corresponding subgroup of the mapping class group.  It must also restrict to the $n=2$ version of $\eta+\lambda\theta$ when restricted to $(n-2)z + \sym^2(\Sigma)$ (where $z$ is a basepoint). Bearing in mind that $H^2(\sym^n(\Sigma);\R)=\R \oplus \Lambda^2 H^1(\Sigma;\R)$ independent of $n\geq 2$, and that $\theta$ and $\eta$ restrict to
$(n-2)z + \sym^2(\Sigma)$ in the obvious way, this shows that it is enough to prove the formula when $n=2$.

Observe that we can embed a one-holed torus $T'$ into $\Sigma$ as a neighbourhood of $\gamma$. The restriction of $\phi^*c_\lambda$ to $\sym^2(T')\subset \sym^2(\Sigma)$ is then independent of the topology of $\Sigma$. From this observation and the case worked out above one obtains the $n=2$ result.
\end{pf}

When $\gamma$ disconnects $\Sigma$, the lemma is still true but rather trivial. In fact, in this case the weight of the blow-up makes no difference to the induced class $c_\lambda$ on $\sym^n(\Sigma)$. These assertions are simple algebraic consequences of the naturality of the construction of $\mathbf{S}_2(\Sigma_\gamma)$ from $\Sigma$ under diffeomorphisms of $\Sigma$ supported away from $\gamma$.

\section{Lagrangian isotopies}\label{isotopies}

We prove statements (1) and (3) of Theorem \ref{main}. The argument also contains one of our two proofs of Proposition \ref{forms}, though this one requires $\lambda>0$ and does not produce K\"ahler forms.

\emph{Step 1a: Reduction to genus 2.}

If $\Sigma$ has genus $g\geq 3$, we can find a simple closed curve $\Gamma$ separating $\Sigma$ into a genus $2$ part containing $\gamma_0$, $\gamma_1$ and $\gamma_2$, and a genus $g-2$ part containing $\gamma_3 \cup \dots \cup \gamma_g$. Thus $\Sigma = \Sigma' \# \Sigma''$, where $\Sigma'$ has genus 2 and contains the handleslide. Fix closed discs $D'\subset \Sigma'$ and $D''\subset \Sigma''$ so that $\Sigma = ( \Sigma' \setminus \interior(D'))\cup_\Gamma (\Sigma'' \setminus \interior(D''))$, and points $p\in \interior(D')$, $q\in \interior(D')$. 

The region of $\sym^g(\Sigma)$ of interest to us is the open set $ \sym^2(\Sigma' \setminus D' ) \times \sym^{g-2}(\Sigma'' \setminus D'')$. We claim that there are K\"ahler forms on $\sym^g(\Sigma)$ which restrict to this region as the sum of forms pulled back from the two factors. Once this claim is established, it will suffice to prove the theorem for a genus 2 surface. However, the reduction will not be perfect: the genus 2 surface has a puncture (far from the handleslide region), and we have to use the symplectic forms produced in the course of proving the claim. 

\begin{lem}
For any $\lambda \geq 0$, there are symplectic (or indeed K\"ahler) forms $\xi$ on $\sym^g(\Sigma' \amalg \Sigma'')$, representing $\eta+\lambda\theta$, such that $i_p^*\xi=i_q^*\xi$. Moreover, $\xi$ restricts to each connected component $\sym^{k}(\Sigma')\times\sym^{g-k}(\Sigma'')$ as the sum of two forms pulled back from the factors. \end{lem}
\begin{pf}
Note that $\eta+\lambda \theta$ is a K\"ahler class (indeed, $\eta$ is ample and $\theta$ is the pullback by the Abel--Jacobi map of an ample class on the Jacobian). Take a K\"ahler form $\kappa'_g$ on $\sym^g(\Sigma')$ representing $\eta+\lambda\theta$ , and let $\kappa'_k$ be its restriction to $\sym^k(\Sigma')$, embedded by $D\mapsto (g-k)p+ D$. Similarly, take a K\"ahler form $\kappa''_g$ on $\sym^g(\Sigma'')$ representing $\eta+\lambda\theta$, and let $\kappa''_k$ be its restriction to $\sym^k(\Sigma'')$, embedded by $D\mapsto (g-k)q+ D$. Define $\xi$, on the connected component $\sym^k(\Sigma')\times \sym^{g-k}(\Sigma'')$, by $\kappa'_{g-k} \oplus \kappa''_k$.
\end{pf}

Given such a form $\xi$, pull it back to $\mathbf{S}_g(\Sigma'\cup \Sigma'')$ (the blow-up of $\sym^g(\Sigma'\cup \Sigma'')$ along $p+q+\sym^{g-2}(\Sigma' \amalg \Sigma'')$) and add a closed form supported near the exceptional divisor so as to obtain a symplectic form $\xi_\lambda$.  The weight of the blow-up is unimportant, because of the remark following the proof of Lemma \ref{coh}. We may assume that $\tilde{\imath}_p^*\xi_\lambda=\tilde{\imath}_q^*\xi_\lambda$, so that $\xi_\lambda$ may be used in forming the symplectic sum $\#_{Z_p \sim Z_q} \mathbf{S}_g(\Sigma' \amalg \Sigma'')$. 

Now, there is a diffeomorphism $f\colon \sym^g(\Sigma) \to \#_{Z_p \sim Z_q} \mathbf{S}_g(\Sigma' \amalg \Sigma'')$ (canonical up to isotopy), and thus $\omega:=f^*\xi_\lambda$ is a symplectic structure on $\sym^n(\Sigma)$. In setting up $f$ we need to fix tubular neighbourhoods of $U_p$ of $Z_p$ and $U_q$ of $Z_q$. It is clear that we can choose $U_p$ and $U_q$ to be disjoint from the open subset
$ \sym^{2}(\Sigma' \setminus D') \times \sym^{g-2}(\Sigma''  \setminus D) 
$ of $\mathbf{S}_g(\Sigma'\amalg \Sigma'')$, and hence we can ensure that $f$ restricts to this subset as the natural inclusion into $\sym^g(\Sigma)$. With this duly arranged, the symplectic form $\omega$ fulfils our claim.

\emph{Step 1b: Completing the reduction to genus 2.} In this step we complete the reduction to genus $2$ by showing that we can choose \emph{any} $j$-positive symplectic form representing the class $\eta+\lambda\theta$ on $\sym^2(\Sigma')$. As things stand, our form is constrained to be the restriction to $(g-2)p + \sym^2(\Sigma')$ of a $j$-positive symplectic form $\kappa_g'$ on $\sym^g(\Sigma')$. 
By an inductive argument, it will suffice to show that we can choose the K\"ahler form freely (within our preferred cohomology class) on the divisor $p+\sym^{g-1}(\Sigma')$. Thus we have a $j$-positive symplectic form on a manifold $M^{2n}$ and  a codimension-2 $j$-holomorphic submanifold $H^{2n-2}$. The boundary of a tubular neighbourhood is a contact type hypersurface (convex, as seen from the inside). But, by an application of Gray's stability theorem---left to the reader---one can always isotope the symplectic form so that its restriction to $H$ is a prescribed $j$-positive symplectic form. This does the trick.

\emph{Step 2: The (punctured) genus 2 case.}
Now assume $\Sigma'$ is closed and has genus 2, and has a basepoint $p$ far from any of the $\gamma_i$. We have to prove the theorem for $\Sigma'\setminus \mathrm{nd}(p)$.

Write $T$ for the torus obtained by surgering out $\gamma_2$. The `scar' left by the surgery is a pair of discs in $T$, containing points $r$ and $s$ in their interiors. In $\Sigma'$, the surgery is done in a region disjoint from the curves $\gamma_0$ and $\gamma_1$, which therefore have images $\bar{\gamma}_0$ and $\bar{\gamma}_1$ in $T$.

At this point we should pin down the choice of $j$-positive symplectic form $\kappa$ on $\sym^2(\Sigma')$. We choose it to be a symplectic structure arising as the pullback via a diffeomorphism $\phi \colon \sym^2(\Sigma') \cong \#_{Z_r \sim Z_s} \mathbf{S}_2(T)$ of a K\"ahler form on $\mathbf{S}_2(T)$. By Lemma \ref{coh}, we can assume that it represents the cohomology class $\eta+\lambda\theta$. 

More specifically, we take a K\"ahler form on $\sym^2(T)$ representing $\eta+\lambda\theta$, and form its ordinary blow up in such a way that the exceptional divisor has weight $\lambda$. We have to make sure that the induced area forms on $Z_r$ and $Z_s$ agree with one another. For this, notice that, by the symplectic neighbourhood theorem, $Z_r$ (say) has a neighbourhood of form $T\times D^2(0;R)$, with a product symplectic form $\varepsilon = \beta+\omega_{C}$. Think of the projection $p \colon T\times D^2(0;R)\to D(0;R)$ as a symplectic fibre bundle. By Thurston's patching argument \cite[Chapter 6]{MS}, we can replace the original symplectic form by a new closed 2-form $ \varepsilon'$, symplectic on the fibres of $p$, agreeing with the $\varepsilon$ over the annulus $R/2 < |z | < R$, such that $\varepsilon'|_{p^{-1}(0)}$ is a freely chosen positive area form $\beta'$. For suitable 2-forms $\zeta$ on $D(0;R)$, supported in $D(0;R/2)$, the form $\varepsilon'+ p^*(\omega_\C+\zeta)$, will then be symplectic. Thus, replacing $\varepsilon$ by $\varepsilon'+ p^*(\omega_\C+\zeta)$ near $Z_r$, we can adjust arrange that the forms on $Z_r$ and $Z_s$ agree.

We also need to be more precise about the symplectic sum and the diffeomorphism $\phi$. To do this, we need to construct suitable tubular neighbourhoods of $Z_r$ and $Z_s$. Given a compact subset $K \subset \Sigma' \setminus \{r,s\} $, there is a natural framing of the restricted normal bundle $N_{Z_r}|_{r+K}$ (more precisely, an identification with the trivial bundle with fibre $T_r \Sigma'$). Thus we can construct a symplectically trivial tubular neighbourhood of $Z_r$ over the subset $r+K\subset Z_r$. Similarly for $Z_s$.

The symplectic sum description provides a hypersurface $V \subset \sym^2(\Sigma')$---the vanishing cycle of the degeneration---whose isotropic foliation is by circles, and whose space of isotropic leaves is identified with $T$.\footnote{The space of isotropic leaves of a (coisotropic) vanishing cycle $V$ is naturally identified with the singular locus in the central fibre of the degeneration. The quotient map to the leaf-space corresponds to the limiting parallel transport map $\rho$ from Proposition \ref{sum}.} Thus there is an $S^1$-bundle  $\rho\colon V\to T$ (topologically trivial, in fact) such that  $\kappa|_V= \rho^*\beta$ for an area form $\beta$ on $T$.  

We want to make sure that $\rho^{-1}(K) = \gamma_2 + K \subset \sym^2(\Sigma)$, and moreover, that $\rho(x+y)=x$ when $x\in K$, $y\in \gamma_2$. 
For this we have to set up the diffeomorphism $\phi$ correctly. We have already set up tubular neighbourhoods of $r+ K$ and $s+K$; points $x+k$ in the boundary of the tubular neighbourhood (with $k\in K$) should map under $\phi$ to $x+k\in \sym^2(\Sigma')$. This is a matter of smooth rather than symplectic topology;  it can be arranged by the method of \cite[Lemma 3.16]{Pe2}, for instance. Note furthermore that, because of the way we have set up the tubular neighbourhoods of $Z_r$ and $Z_s$, we end up with a symplectic form $\kappa$ on $\sym^2(\Sigma)$ which is product-like in $K\times \mathrm{nd}(\gamma_2)$.
 
\begin{rmk}
Our symplectic forms on $\sym^g(\Sigma)$ are not necessarily K\"ahler; the symplectic sum operation does not in general preserve the K\"ahler category. However, it is easy to see that they can be taken to tame the natural complex structure $\sym^n(j)$. 
\end{rmk}

\emph{Step 3: Lagrangian isotopy.}
Take an isotopy $\{\bar{\gamma}_t\}_{t\in [0,1]}$ of circles in $T$, from $\bar{\gamma}_0$ to $\bar{\gamma_1}$. It is, of course, a Lagrangian isotopy with respect to $\beta$. The preimages $\rho^{-1}(\bar{\gamma}_t)$ are tori  in $\sym^2(\Sigma')$. They are Lagrangian, because if $u$ and $v$ are tangent to $\rho^{-1}(\bar{\gamma}_t)$ at a point $\rho^{-1}(z)$ then $\kappa(u,v) = \beta(\rho_*u,\rho_* v ) =0$. Moreover, we set things up in Step 2 so that $\rho^{-1}(\bar{\gamma_i})= \gamma_i \times \gamma_2$ for $i=0,1$. Thus we have a Lagrangian isotopy from $\gamma_0 \times \gamma_2$ to $\gamma_1\times \gamma_2$. Notice that this isotopy stays inside a compact subset of $\sym^2(\Sigma' \setminus \{p\})$, so it is again compatible with step 1.
 
This completes the proof of clauses (1) and (3) of Theorem \ref{main}.

\emph{Step 4: Hamiltonian isotopy.} We begin the proof of statement (2) of Theorem \ref{main}. 

The difference between Lagrangian and Hamiltonian isotopy is measured by the \emph{flux} \cite{MS}. Given a Lagrangian isotopy $\{L_t\}$, its (infinitesimal) flux at time $t$ is a class $a_t\in H^1(L_t;\R)$ which can be obtained by extending the vector field generating the isotopy to a globally-defined symplectic vector field, dualising this to get a closed 1-form, and taking its cohomology class restricted to $L_t$. For the isotopy in Step 3, one checks that the flux $a_t$ is equal to $\rho^*b_t$, where $b_t$ is the flux of the isotopy $\{\bar{\gamma}_t\}$ in $T$ with respect to $\beta$. Thus our Lagrangian isotopy in $\sym^2(\Sigma')$ is Hamiltonian if and only if $\{\bar{\gamma}_t\}$ is  a Hamiltonian isotopy in $T$.  Now,
the complement $T \setminus (\bar{\gamma}_0\cup \bar{\gamma}_1)$ has three components, of which two are homeomorphic to the open disc. The only obstruction to making $\{\bar{\gamma}_t\}$ Hamiltonian is that the two disc-components have equal $\beta$-area.

We have now proved the existence of Hamiltonian isotopies under the condition that two discs have the same area with respect to the area form $\beta$ on the surgered surface. However, the stated criterion concerned areas with respect to the area form $\alpha$ on $\Sigma$, and we need to relate these conditions. In stating the theorem, we noted that the existence of a Hamiltonian isotopy depends only on $\lambda$ and $\alpha$ (and, \emph{a priori}, on $j$). The problem has an extra symmetry, however: we can apply self-diffeomorphisms to $\Sigma$, pulling back $j$,  provided that they act trivially on $\gamma_0\cup \gamma_1 \cup \gamma_2$. 
\begin{lem}
The only invariants of area forms $\alpha$ under such diffeomorphisms are the areas of $A$, $D$ and $\Sigma \setminus (A \cup D)$.
\end{lem}
\begin{pf}
Given area-forms $\alpha_0$ and $\alpha_1$, we can certainly find a self-diffeomorphism $\psi$, trivial along $\Gamma:= \gamma_0\cup \gamma_1\cup \gamma_2$, such that $\psi^*\alpha_0 = \alpha_1$ in a neighbourhood of $\Gamma$. Now replace $\alpha_0$ by $\psi^* \alpha_0$ and apply Moser's deformation argument to the path $\alpha_t = (1-t)\alpha_0+ t\alpha_1$.
\end{pf}
We thus see that the existence of Hamiltonian isotopies depends on $\alpha$ only through the areas of $A$ and $D$ (the area of their complement is clearly irrelevant). Our construction can be interpreted, then, as saying that there is a function $f(\lambda,\int_D \alpha)$ such that a Hamiltonian isotopy exists if $\int_A\alpha = f(\lambda,\int_D\alpha)$. In Section \ref{non-iso} we shall show that a Hamiltonian isotopy can \emph{only} exist if $\int_A\alpha = (1+\lambda)^{-1}(\int_D \alpha+ 2\lambda) $, which must therefore coincide with $ f(\lambda,\int_D\alpha)$.

\section{Hamiltonian non-isotopy}\label{non-iso}
In this section we complete the proof of Theorem \ref{main} by showing that the tori $T_0$ and $T_1$ can only be Hamiltonian-isotopic when $\int_D{\alpha}$ is related in a precise way to $\int_A{\alpha}$. 

We do this by showing that if the area constraint fails there is a discrepancy between the ranks of the Lagrangian Floer homology groups $HF_*(T_0,T_0)$ and $HF_*(T_0,T_1)$. These would be isomorphic if a Hamiltonian isotopy existed. We  work over the universal Novikov ring $\Lambda_{\Z/2}$ of the field $\Z/2$, i.e., the ring of formal `series' $\sum_{r\in \R}{a(r)t^r}$ where $a\colon \R\to \Z/2$ is a function such that $\mathrm{supp}(a)\cap (-\infty, c]$ is finite for all $c$. This coefficient ring is used to record the areas of holomorphic discs. 
The well-definedness and invariance of Floer homology for monotone Lagrangians of minimal Maslov index 2 is not automatic, but in the case of Heegaard tori it follows from a cancellation theorem from \cite{OS} (Theorem 3.15) alongside the transversality argument outlined in \cite{Oh}. 

The self-Floer homology of $T_0$ is as large as it conceivably could be:
\begin{lem}
 $HF_*(T_0,T_0)\cong H_*( T_0 ;\Lambda_{\Z/2}) $ for any form $\omega_\lambda$ as in Proposition \ref{forms}.
\end{lem}
\begin{pf}
The relevant moduli spaces are computed by Ozsv\'ath--Szab\'o in \cite[Lemma 9.1]{OS}. The discs contributing to the Floer-theoretic differential come in pairs $(u_1,u_2)$, with $u_1$ and $u_2$ elements of the same moduli space. These come from discs in $\Sigma$ itself, of equal $\alpha$-areas and disjoint from the diagonal $\delta\subset \sym^g(\Sigma)$. Thus they have equal $\omega_0$-area when $\omega_0$ is a form as in Proposition \ref{forms} representing $\eta$. When $\lambda = 0$ the lemma therefore follows from Ozsv\'ath--Szab\'o's calculation.

Now perturb $\omega_0$ to $\omega_\lambda$ with $[\omega_\lambda]=\eta+\lambda\theta$ by adding a small multiple of a closed $(1,1)$-form supported near $\delta$ and representing the dual of $[\delta]$. This perturbation does not change the areas, and so $u_1$ still cancels $u_2$ and the argument goes through.
\end{pf}
The Floer homology of $(T_0,T_1)$ is usually smaller:
\begin{lem}
Either $\rank_{\Lambda_{\Z/2}} HF_*(T_0, T_1) =2^{g-2}$, or
$\int_D{\alpha} = (1+\lambda) \int_A{\alpha} -2 \lambda.$ 
\end{lem}
\begin{pf}
Again, the proof rests on a calculation from \cite{OS} (Lemma 9.4). We set things up as in \cite[Figure 9]{OS}, intersecting $\gamma_0\times \gamma_2$ with $\gamma_1\times \gamma_2'$, where  $\gamma_2'$ is a small, transverse, Hamiltonian perturbation of $\gamma_2$ such that $\gamma_2\cap \gamma_2'=2$. Let $\gamma_2\cap \gamma_2'= \{ u,v \}$ and $\gamma_0\cap \gamma_1 = \{x,y\}$. Floer's complex is then a tensor product $C' \otimes_{\Lambda_{\Z/2}} C''$, where $C''$ (the `uninteresting' part) corresponds to intersections among the curves $(\gamma_3,\dots ,\gamma_g)$ and their $\alpha$-Hamiltonian translates, while $C'$ has four generators: $\mathbf{x} = u+ x$, $\mathbf{y}=u+y$, $\mathbf{x}'=v+ x$ and $\mathbf{y}'=v+y$. We set up the labelling so that the only potentially non-zero matrix entries in the differential $\partial$ on $C'$ are $\langle \partial \mathbf{x}, \mathbf{y} \rangle$ and  $\langle \partial \mathbf{x}', \mathbf{y}'\rangle$. Now, writing elements of $\Lambda_{\Z/2}$ as sums $\sum_{r\in \R}{n_r t^r}$ with $n_r\in \Z/2$, Ozsv\'ath--Szab\'o's calculation shows that 
\[\langle \partial \mathbf{x}, \mathbf{y} \rangle = 
 t^{\int_{D_1}\omega_\lambda} + t^{\int_{D_2}\omega_\lambda}\] 
where $D_1$ and $D_2$ are two particular discs and the integrals are with respect to $\omega_\lambda$. The first disc, $D_1$, is the product of the disc $D$ shown in Figure 1 and a constant disc at $u$. It does not intersect the diagonal $\delta$, and so (using perturbed K\"ahler forms $\omega_\lambda$ as in the proof of the last lemma) its area is $\int_D\alpha$, regardless of $\lambda$. The second disc, $D_2$, arises from a branched double covering of a disc by the annulus $A$. By Riemann--Hurwitz, its intersection number with $\delta$ is 2. Hence, using the formula $[\delta] =  \eta  -  \theta$ (which was mentioned following Proposition \ref{forms}), we have
\[ \int_{D_2}{\omega_\lambda} =   \int_A{\alpha} + \lambda \left( \int_A{\alpha} -2 \right) .  \] 
The homology of $C'$ is zero unless the areas of $D_1$ and $D_2$ are equal, that is, unless $ \int_{D}\alpha  = (1+\lambda) \int_A{\alpha}  -2\lambda$. But $H_*(C'')$ has rank $2^{g-2}$ (as in the proof of the previous lemma), so the result follows.
\end{pf}
The `only if' clause in statement (2) of Theorem \ref{main} follows immediately. The `if' clause also follows, by the argument explained in the last paragraph of Section \ref{isotopies}.

\section{Application to Heegaard Floer homology}\label{Heegaard}
We now turn to the proof of Corollary \ref{HeegCor}, freely drawing on the language of Heegaard Floer homology developed by Ozsv\'ath--Szab\'o in \cite{OS}. 
Before we begin the proof, we make some remarks on the relation of Heegaard Floer theory to Lagrangian Floer homology. 

It is a feature of Lagrangian Floer homology that the moduli spaces defining the differential depend on the Lagrangians and the almost complex structure but do not directly involve a symplectic form. The primary role of symplectic forms is to fix the area of the pseudo-holomorphic `Whitney discs' in a given homotopy class. Having an \emph{a priori} bound on area is essential for compactness, and for this reason it is important to work with Lagrangian (not merely totally real) submanifolds. Lacking a convenient symplectic form, Ozsv\'ath and Szab\'o bound the areas of Whitney discs by treating $\sym^g(\Sigma)$ as a `symplectic orbifold'---the quotient of $\Sigma^{\times g}$, with its product symplectic form, by the action of the symmetric group---and estimating areas `upstairs'. 

The symplectic class plays a second role in Floer theory: the areas of the holomorphic discs are typically recorded in a Novikov ring of coefficients, just as in Section \ref{non-iso}. This is useful because there may be infinitely many homotopy classes of Whitney discs, but there are only finitely many index 1  Whitney discs with a given area.

Ozsv\'ath--Szab\'o are able to do without Novikov rings by showing that, for a fixed $\spinc$-structure $\mathfrak{s}$, and under the assumption of `weak admissibility' of the pointed Heegaard diagram $(\Sigma,\boldsymbol{\alpha},\boldsymbol{\beta},z)$ with respect to $\mathfrak{s}$, there are only finitely many homotopy classes of Whitney discs, between intersection points $\mathbf{x}\in \mathbb{T}_\alpha \cap \mathbb{T}_\beta$ with $s_z(\mathbf{x})=\mathfrak{s}$, with a given Maslov index, fixed intersection number with the divisor $z + \sym^{g-1}(\Sigma)$, and additionally satisfying a positivity constraint automatically satisfied by pseudo-holomorphic discs \cite[Lemma 4.14]{OS}.

Now suppose $\omega_\lambda$ is as in Proposition \ref{forms}. Lagrangian Floer homology $HF(\mathbb{T}_\alpha, \mathbb{T}_\beta)$, in any of its possible algebraic variants, can be computed via any almost complex structure taming $\omega_\lambda$ that satisfies a regularity condition. Ozsv\'ath and Szab\'o show that to attain regularity it is sufficient to consider (paths in) a particular class of almost complex structures, those that are small `nearly symmetric' perturbations of the standard integrable structure \cite[Theorem 3.4]{OS} (these still tame $\omega_\lambda$). The periods are controlled by the Maslov index and intersection numbers with $z+\sym^{g-1}(\Sigma)$, and the groups are set up so as to keep track of these intersection numbers.

Heegaard tori are monotone Lagrangians of minimal Maslov index 2. The troublesome feature of Lagrangian Floer theory in such a case is that bubbling-off of discs spoils compactness of the moduli spaces involved in the definition of the groups (more precisely, the matrix coefficients $\langle \partial \circ \partial x , x\rangle$ where $x$ is a generator for the Floer complex and $\partial$ the usual Floer-theoretic `differential'), cf. \cite{Oh}. Ozsv\'ath and Szab\'o show that, though Maslov index 2 discs can bubble off in the relevant 1-dimensional moduli space, they always do so in cancelling pairs \cite[Theorem 3.15]{OS}, which, as explained  in \cite{Oh}, implies that $\partial \circ \partial =0$.

\begin{pf}[Proof of Corollary \ref{HeegCor}]
We shall set up a  `continuation isomorphism' 
\[HF^+( \boldsymbol{\alpha}, \boldsymbol{\beta },\mathfrak{s})\to HF^+(
\boldsymbol{\alpha}, \boldsymbol{\beta },\mathfrak{s}),\] 
and a similar one for the other variants of Heegaard Floer homology.

By hypothesis $\mathbb{T}_{\alpha}$ is transverse to $\mathbb{T}_{\beta}$, and perturbing the $\boldsymbol{\alpha}'$ curves slightly we may assume that $\mathbb{T}_{\alpha}$ is transverse to $\mathbb{T}_{\beta}$. Choose any small $\lambda>0$ and an area form $\alpha$ satisfying the area constraint of Theorem \ref{main} (2). We then consider a symplectic form $\omega_\lambda$ as in Proposition \ref{forms} and use Theorem \ref{main} to obtain a Hamiltonian isotopy $\{  \phi_t\}_{t\in [0,1]} $ with $\phi_0=\id$ and $\phi_1(\mathbb{T}_\alpha)=\mathbb{T}_{\alpha'}$. 

By a general principle in Lagrangian Floer homology, the Hamiltonian isotopy $\{\phi_t\}_{t\in [0,1]}$ induces a continuation isomorphism $HF_*(\mathbb{T}_\alpha,\mathbb{T}_\beta;\omega_\lambda) \to HF_*(\phi_1(\mathbb{T}_\alpha), \mathbb{T}_\beta;\omega_\lambda)$ (Floer homology with Novikov-ring coefficients). Because the minimal Maslov index is 2, the continuation-map moduli spaces cannot bubble here (and the same goes for those used in proving the composition law for continuation maps).  However, it is not immediately apparent that the continuation map makes sense for $HF^+$ (say), because the $\omega_\lambda$-area of Whitney discs might not be controlled by Maslov index and intersection number with $z+\sym^{g-1}(\Sigma)$. We deal with this issue now.

The continuation map can be understood as follows.  We have a Hamiltonian isotopy $\{\phi_s\}_{s\in [0,1]}$, and we may extend this to a family $\{\phi_s\}_{s\in \R}$ where $\phi_s = \id$ for $s<0$ and $\phi_s=\phi_1$ for $s>1$ (this family will not be smooth, so strictly we should reparametrise $[0,1]$ before we extend so as to obtain a smooth family). We consider the trivial $\sym^g(\Sigma)$-bundle $\sym^g(\Sigma)\times \R\times [0,1] \to  \R\times [0,1]$, and make it into a Hamiltonian fibration by giving endowing it with a closed 2-form $\Omega:=\omega_\lambda+ d(\beta(t)H_s ds)$. Here $s$ is the $\R$-coordinate, $t$ the $[0,1]$-coordinate; $\{H_s\}_{s\in \R}$ is a family of functions, non-zero only for $s\in [0,1]$, generating $\phi_s$; and $\beta(t)$ is a cut-off function, equal to $1$ near $t=1$ and to $0$ near $t=0$. There are Lagrangian boundary conditions $\mathbb{T}_{\beta}\times \R\times \{1\}$ and $\bigcup_{s\in \R}{\phi_s( \mathbb{T}_\alpha) \times \{ s\} \times \{0 \}}$ over the boundary of the strip. The continuation map is defined by counting index-0 pseudo-holomorphic sections of this fibration subject to the Lagrangian boundary conditions (for an almost complex structure making the projection holomorphic, and translation-invariant for $|s|\gg 0$).  See \cite{Sei2} for details of the fibre-bundle approach to Floer theory.

We need to estimate the energies of sections $u$ subject to the Lagrangian boundary conditions and asymptotic to intersection points $\mathbf{x}_-\in \mathbb{T}_\alpha \cap \mathbb{T}_{\beta}$ as $s\to -\infty$ and $\mathbf{x}_+\in  \phi_1(\mathbb{T}_\alpha) \cap \mathbb{T}_{\beta}$ as $s\to +\infty$. 

These energies are cohomological in nature; they are invariant under compactly supported homotopies of $u$, for instance. Because of this, we can estimate the energies by considering a degeneration of the domain in which a disc is pinched off, as in Figure 2. That is, we consider a 1-parameter family of surfaces $\{S_t\}_{t\in [0,\infty)}$, with $S_0=\R\times [0,1]$, pulling back our Hamiltonian fibration and Lagrangian boundary conditions by a family of diffeomorphisms $f_t\colon S_t\to S_0$. We arrange that the surfaces $S_t$ limit to the nodal surface $S_\infty$ shown in Figure 2 as $t\to \infty$. Conformally, we can view $S_\infty$ as the nodal union of a triangle $T$ and a disc $D$. We can arrange that there is a limiting Hamiltonian fibration over $S_\infty$ (still smoothly trivialised), and that over the three edges of the triangle, we have constant Lagrangian boundary conditions $\mathbb{T}_\alpha$, $\phi_1(\mathbb{T}_\alpha)$ and $\mathbb{T}_\beta$, as shown. Over the boundary of the disc we have the Lagrangian boundary condition $\phi_s(\mathbb{T}_\alpha)$, where $s$ now parametrises the boundary anticlockwise.

Consider triangles $T\to \sym^g(\Sigma)$, subject to the boundary conditions given by  $\mathbb{T}_\alpha$, $\mathbb{T}_{\alpha'}=\phi_1(\mathbb{T}_\alpha)$ and $\mathbb{T}_\beta$, of Maslov index $0$, and with fixed intersection number with $z+\sym^{g-1}(\Sigma)$.  Assume that $T$ is asymptotic to intersection points $\mathbf{x}_\pm$, as above, whose associated $\spinc$-structure is $\mathfrak{s}$. According to \cite[Section 9.2]{OS}, the number of finite-energy holomorphic triangles $T$ satisfying these conditions is finite. On the other hand, homotopy classes of sections over $D$, subject to the boundary conditions $\phi_s(\mathbb{T}_\alpha)$ and mapping the node to a fixed $\mathbf{x}_0\in \mathbb{T}_\alpha \cap \phi_1(\mathbb{T}_\alpha)$, correspond to $\pi_2(\sym^g(\Sigma), \mathbb{T}_\alpha)=\Z$ (by extending them to sections of a trivial fibration over a larger disc, with constant Lagrangian boundary condition $\mathbb{T}_\alpha$), and so again the energy is controlled by the intersection number with $z+\sym^{g-1}(\Sigma)$. 

Each homotopy class of sections over $S_\infty$ can be smoothed out to a  homotopy class of sections of $S_t$ for $t$ finite, and it is easy to see that every homotopy class of sections over $S_0$ arises this way. The smoothing leaves the energy unchanged. Hence, going back to the continuation map picture (over $S_0$) we conclude that the energy of index-0 pseudo-holomorphic strips with fixed asymptotics and fixed intersection number with $z+\sym^{g-1}(\Sigma)$ is bounded. Thus the continuation map is well-defined. More accurately, the continuation map is defined initially on $CF^\infty$, but it respects the subcomplexes $CF^-$ and so induces maps on $CF^-$ and $CF^+$. For $\widehat{CF}$ one defines the continuation map using strips whose intersection number with $z+\sym^{g-1}(\Sigma)$ is zero.

By using the same degeneration to obtain the necessary energy bounds, one sees that these continuation maps are chain maps, and that the continuation map associated with the reverse isotopy $\phi_{1-t}$ is an inverse up to chain-homotopy. Hence all the continuation maps are quasi-isomorphisms.  

Note finally that, because Theorem \ref{main} produces Hamiltonian isotopies which are canonical (up to homotopy with fixed endpoints), the continuation isomorphisms on Heegaard Floer homology are also canonical.
\end{pf}
\begin{figure}[t!]
\labellist
\small\hair 2pt
\pinlabel $\mathbb{T}_\beta$ at 280 605
\pinlabel $\mathbb{T}_\alpha$ at 200 515
\pinlabel $\phi_1(\mathbb{T}_\alpha)$ at 370 515
\pinlabel $\phi_s(\mathbb{T}_\alpha)$ at 215 450
\pinlabel $T$ at 280 560
\pinlabel $D$ at 280 460
\endlabellist
\centering
\includegraphics[width=0.6\textwidth]{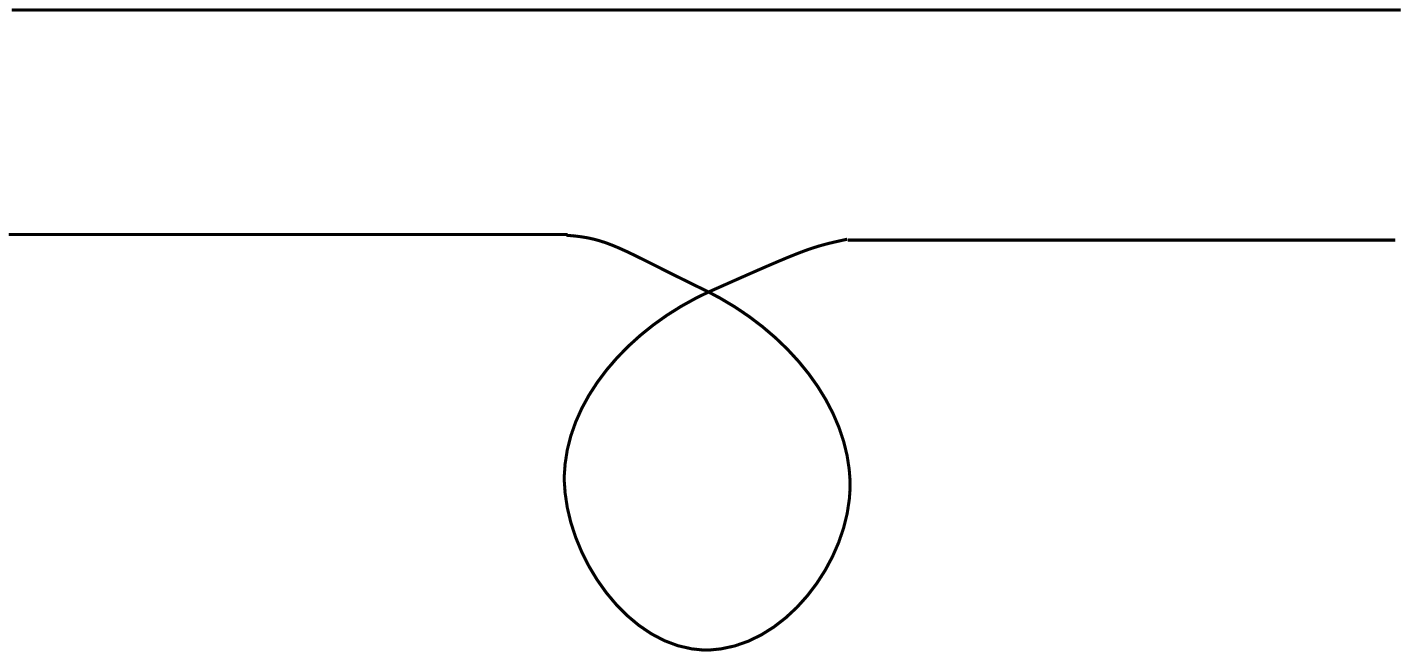}
\caption{Degeneration of the domain to the union of a (conformal) triangle and a disc.}
\end{figure}

\begin{rmk}
The continuation isomorphisms on Heegaard Floer homology which we have associated with a handleslide are the same as the isomorphisms constructed by Ozsv\'ath--Szab\'o using holomorphic triangles. To see this one applies a gluing theorem for holomorphic sections to the degeneration of Figure 2.
\end{rmk}

\section{Branched coverings and smoothing of currents}\label{smoothing}

This final section reproduces the note \cite{Pe1}, which has not previously been published. 

We consider branched coverings $\pi\colon X\to X' $ of complex manifolds---that is, holomorphic maps which are proper, surjective, and finite. The branch locus $B_\pi \subset X'$ of such a map is 
\[ B_\pi = \{ \pi(x): x\in X, \; \ker D_x(\pi) \neq 0\}.\] 
A $C^\infty$ K\"ahler form $\omega$ on $X$ can be pushed forward---in the sense of currents, that is, of 2-forms with coefficients of class  $L^1_{\mathrm{loc}}$---to a closed current $\pi_*(\omega)$ on $X'$ which is smooth on $X'\setminus B_\pi$. 
The following result is essentially due to Varouchas \cite{Var}.
\begin{prop}\label{var}
Let $\pi\colon  X \to X'$ be a branched covering of complex manifolds,
and $\omega$ a K\"ahler form on $X$. Let $N$ be a
neighbourhood of the branch locus in $X'$. Then there exists a
K\"ahler form $\omega'$ on $X'$, representing the class $[\omega']=\pi_*[\omega]\in H^2(X';\R)$, such that $(\pi_*\omega-\omega')|_{X'\setminus \overline{N}}=0$. 
\end{prop}
We explain here the minor modification of Varouchas' argument from \cite{Var} needed to prove this result (the stated conclusion in \cite{Var} is simply that $X'$ admits a K\"ahler form). 

Proposition \ref{var} has the following special case:
\begin{cor}\label{smoothing cor}
Let $\Sigma$ be a Riemann surface with positive area form $\alpha$. Let 
$ \pi \colon \Sigma^{\times n}\to \sym^n(\Sigma)$ be the projection map. Suppose that $N\subset \sym^n(\Sigma)$ is an open subset containing the  diagonal $\delta$. Then there exists a K\"ahler form $\omega$ on $\sym^n(\Sigma)$, representing the class $\eta=\pi_*[\alpha^{\times n}]$, such that
outside $\overline{N}$, $\omega$ is the smooth push-forward $\pi_*(\alpha^{\times n})$ of the product form.
\end{cor}
\begin{defn}
Let $X$ be a complex manifold. A {\bf K\"ahler cocycle} on $X$ is
a collection $(U_i,\varphi_i)_{i\in I}$, where $(U_i)_{i\in
I}$ is an open cover of $X$, and $\varphi_i\colon U_i \to \R$ is a
a function, such that for all $i,j \in I$,
\begin{enumerate}
\item $\varphi_i$ is strictly plurisubharmonic on $U_i$; and
\item $\varphi_i-\varphi_j$ is pluriharmonic on $U_i \cap U_j$.
\end{enumerate}
\end{defn}
One ascribes to the cocycle a property (continuity, smoothness,
etc.) possessed by all the $\varphi_i$. K\"ahler cocycles are, by
definition, upper semicontinuous.
Condition (1) means that the $2$-current $dd^c \varphi_i$ is
strictly positive on $U_i$; (2) means that these currents agree on
overlaps, and are therefore restrictions of a $2$-current $\omega$
on $X$ (closed and strictly positive). If the cocycle is
$C^\infty$ then $\omega$ will be a K\"ahler form.

Varouchas' \emph{lemme principal} is the following. The proof uses the `regularised maximum' technique of Richberg and Demailly.
\begin{lem}
Let $U,V,W,\Omega$ be open subsets of $\C^n$ with
\[U\Subset V\Subset W, \quad \Omega \subset W. \]
Let $\phi\colon W \to \R$ be continuous, strictly plurisubharmonic, and smooth on $\Omega$. 
Then there exists a function $\psi \colon W\to \R$, again continuous and strictly plurisubharmonic, equal to $\phi$ on $W \setminus \overline{V}$ and smooth on $U\cup \Omega$.
\end{lem}
(The notation $U\Subset V$ means that $\overline{U}\subset V$.) One passes from local to global by the following argument, which I give in detail since Varouchas' stated conclusion is weaker. 
\begin{lem}\label{global smoothing}
Let $(U_i,\varphi_i)_{i\in I}$ be a continuous K\"ahler cocycle on the complex manifold $X$. Suppose that $X=X_1\cup X_2$, with $X_1$ and $X_2$ open,
and that the functions $\varphi_i|_{U_i \cap X_1}$ are smooth. Then
there exists a continuous function $\chi \colon X \to \R$ supported in $ X_2$
and a locally finite refinement $ V_j \subset U_{i(j)} \quad (j \in J) $ so that the family
\[   \left(V_j, \varphi_{i(j)}|_{V_j}+\chi|_{V_j} \right)_{j\in J}     \]
is a smooth K\"ahler cocycle.
\end{lem}
\begin{pf}
Refine the cover $(U_i)_{i\in I}$ to a countable, locally finite
cover $(V_i)_{i\in I_1 \amalg I_2}$ with the property that
\[ i \in I_\alpha \;  \Rightarrow  \;V_i \subset X_\alpha, \quad \alpha=1,2. \]
For definiteness suppose both $I_1$ and $I_2$ are infinite;
say $I_\alpha = \mathbb{N}\times \{\alpha\}$, so that the labels $i$ are pairs $(k,1)$ with $k\in \Z$ for $I_1$, or $(k,2)$ for $I_2$. Find open subsets $V_i'' \Subset V_i' \Subset V_i $ such that $(V_i'')$ still covers $X$. Set
\[  A_1=\emptyset, \quad A_n = V''_{(1,2)}\cup\dots\cup V''_{(n-1,2)} \quad (n>1). \]
The sets $A_n$ exhaust $X_2\setminus X_1$. Let
$(V_i,\psi^1_i)_{i\in I_1\cup I_2}$ be the K\"ahler cocycle
induced from $(U_i,\varphi_i)_{i\in I}$ by the
refinement.

\emph{Claim:} there are K\"ahler cocycles $(V_i,\psi^n_i)$, where $n=1,2,\dots$ indexes the elements of
$I_2$, such that the following hold for all $i \in I_1 \cup I_2$ and all $n>1$:
\begin{enumerate}
\item   $\psi^n_i$ is smooth on the set $V_i \cap (X_1\cup A_n).$
\item   There is a continuous function $\chi_n\colon X \to \R$, with $\mathrm{Supp}(\chi_n)\subset
        V'_{(n-1,2)}$, such that $\psi^n_i=  \psi^{n-1}_i +\chi_n$.
\end{enumerate}
We prove the claim by induction on $n$. Apply the `lemme principal' to
\[ (U,V,W,\Omega) = (V''_{(n-1,2)},V'_{(n-1,2)},V_{(n-1,2)}, V_{(n-1,2)}\cap(X_1\cup A_{n-1})) \]
and to the function $\psi^{n-1}_{(n-1,2)}$, obtaining a new function $\psi^{n}_{(n-1,2)} $; let $\chi_n = \psi^{n}_{(n-1,2)} - \psi^{n-1}_{(n-1,2)}$, extended by zero to all of $X$, and use (2) to define the new cocycle. We have to verify (1), i.e. to prove smoothness of $\psi^n_i$ at each $x\in V_i\cap(X_1 \cup A_{n})$. If $x\not \in V'_{(n-1,2)}$ then $\chi_n(x)=0$, but $\psi_i^{n-1}$ was already smooth. If $x\in V_{(n-1,2)}'$ then, near $x$, $\psi^n_i = (\psi^n_i -\psi^n_{(n-1,2)})+ \psi^n_{(n-1,2)} = (\psi^1_i - \psi^1_{(n-1,2)}) + \psi^n_{(n-1,2)}$, which is the sum of a pluriharmonic function and a smooth plurisubharmonic one. But a pluriharmonic function is automatically smooth. By a similar argument, $\psi^n_i$ is strictly plurisubharmonic.

Now define a function $\chi\colon X \to \R$ by $ \chi(x) = \sum_{n \geq 1}{\chi_n(x)}$ (the sum is locally finite).
Then $ \psi^\infty_i(x)  :=  \psi^1_i(x) +\chi(x)$
defines a K\"ahler cocycle. It is
smooth, since on one hand, on $X \setminus {\bigcup{V'_{(n,2)}}} \subset X_1$, the original cocycle was smooth and has not been modified, while on the other hand,
\[ V_{(n,2)} \subset X_1 \cup \bigcup{A_k},  \]
so smoothness on $V_{(n,2)}$ is guaranteed by (1). Hence $\chi$ has the
required properties.
\end{pf}
\begin{pf}[Proof of Proposition \ref{var}]
Each fibre $\pi^{-1}(x')$, being finite, has a neighbourhood which is a disjoint union of open balls. Hence, using the $dd^c$-lemma, one can find a smooth K\"ahler cocycle $(U_i,\varphi_i)$ on $X$ such that each $U_i$ contains a fibre of $\pi$, with $\omega|_{U_i}=dd^c \varphi_i$. 
One can then find a locally finite cover $(U'_i)$ of $X'$ such that
$U_i \supset \pi^{-1}(U'_i)$.

A general property of branched covers is that the push-forward $\pi_*f $ of a continuous function $f\colon X\to \R$ is again continuous (it is given by $\pi_*f (x')=\sum_{x\in\pi^{-1}(x')}f(x)$, where the points $x$ are taken with multiplicities). The family $(U'_i,\pi_*\varphi_i)$ on $X'$ is thus a continuous K\"ahler cocycle: plurisubharmonicity is clear away from $B_\pi$, since $\pi_* dd^c f = dd^c\pi_*f $, hence everywhere by density; similarly for pluriharmonicity on overlaps. As for strictness, let $\omega_{\C^n}$ be the standard K\"ahler form in a complex chart centred at a point $x\in B_\pi$. Then, for a test $(2n-2)$-form $\beta$  supported near $x$, one has $(\pi_*dd^c f - \epsilon \omega)(\beta) =  \pi_*(dd^c f -  (\deg \pi)^{-1}\epsilon\pi^* \omega_{\C^n})(\eta) = (dd^c f-(\deg\pi)^{-1}\epsilon\pi^* \omega_{\C^n})(\pi^*\eta)$, and using this one verifies strict positivity of $\pi_*dd^c f$.

Now let $N' \Subset N$ be a smaller open neighbourhood of $B_\pi$.
Apply the `global smoothing' lemma (\ref{global smoothing}) to $(U'_i,\pi_*\varphi_i)$ on
$X'$, taking $X_1=X'\setminus \overline{N'}$ and $X_2=N$. The
output is a function $\chi\colon X' \to \R$, as well as a
refinement of $(U_i')$, such that
\[  \omega_{X'}:= dd^c (\pi_*\varphi_i +\chi)= \pi_*dd^c\varphi_i+dd^c\chi \]
is a well-defined $2$-form with the right properties. For any smooth, closed test form $\beta\in \Omega^{2n-2}_c(X')$ one has  $\int_{X'} dd^c \chi \wedge \beta =0 $, hence $dd^c\chi$ represents the zero cohomology class.
\end{pf}

\end{document}